\documentclass[twoside,leqno,11pt, A4]{amsart}
\usepackage{amsfonts}
\usepackage{amsmath}
\usepackage{amscd}
\usepackage{amssymb}
\usepackage{amsthm}
\usepackage{amsrefs}
\usepackage{latexsym}
\usepackage{mathrsfs}
\usepackage{bbm}
\usepackage{enumerate}
\usepackage{color}
\begin{document}

\newtheorem{theorem}[subsection]{Theorem}
\newtheorem{proposition}[subsection]{Proposition}
\newtheorem{lemma}[subsection]{Lemma}
\newtheorem{corollary}[subsection]{Corollary}
\newtheorem{conjecture}[subsection]{Conjecture}
\newtheorem{prop}[subsection]{Proposition}
\numberwithin{equation}{section}
\newcommand{\mr}{\ensuremath{\mathbb R}}
\newcommand{\mc}{\ensuremath{\mathbb C}}
\newcommand{\dif}{\mathrm{d}}
\newcommand{\intz}{\mathbb{Z}}
\newcommand{\ratq}{\mathbb{Q}}
\newcommand{\natn}{\mathbb{N}}
\newcommand{\comc}{\mathbb{C}}
\newcommand{\rear}{\mathbb{R}}
\newcommand{\prip}{\mathbb{P}}
\newcommand{\uph}{\mathbb{H}}
\newcommand{\fief}{\mathbb{F}}
\newcommand{\majorarc}{\mathfrak{M}}
\newcommand{\minorarc}{\mathfrak{m}}
\newcommand{\sings}{\mathfrak{S}}
\newcommand{\fA}{\ensuremath{\mathfrak A}}
\newcommand{\mn}{\ensuremath{\mathbb N}}
\newcommand{\mq}{\ensuremath{\mathbb Q}}
\newcommand{\half}{\tfrac{1}{2}}
\newcommand{\f}{f\times \chi}
\newcommand{\summ}{\mathop{{\sum}^{\star}}}
\newcommand{\chiq}{\chi \bmod q}
\newcommand{\chidb}{\chi \bmod db}
\newcommand{\chid}{\chi \bmod d}
\newcommand{\sym}{\text{sym}^2}
\newcommand{\hhalf}{\tfrac{1}{2}}
\newcommand{\sumstar}{\sideset{}{^*}\sum}
\newcommand{\sumprime}{\sideset{}{'}\sum}
\newcommand{\sumprimeprime}{\sideset{}{''}\sum}
\newcommand{\shortmod}{\ensuremath{\negthickspace \negthickspace \negthickspace \pmod}}
\newcommand{\V}{V\left(\frac{nm}{q^2}\right)}
\newcommand{\sumi}{\mathop{{\sum}^{\dagger}}}
\newcommand{\mz}{\ensuremath{\mathbb Z}}
\newcommand{\leg}[2]{\left(\frac{#1}{#2}\right)}
\newcommand{\muK}{\mu_{\omega}}
\newcommand{\seteq}{:=}
\newcommand{\odd}{\mathrm{\ primary}}
\newcommand{\res}{\mathrm{Res}}

\title[Large values of quadratic Dirichlet $L$-functions of prime-related moduli]{Large values of quadratic Dirichlet $L$-functions of prime-related moduli}

\author[P. Gao]{Peng Gao}
\address{School of Mathematical Sciences, Beihang University, Beijing 100191, China}
\email{penggao@buaa.edu.cn}

\begin{abstract}
In this paper, we apply the long
resonator method to exhibit large values at the central point of the family of quadratic Dirichlet $L$-functions of prime-related moduli under the generalized Riemann hypothesis. 
\end{abstract}

\maketitle



\section{Introduction}

Much progress has been made towards understandings of maximum sizes of $L$-functions ever since the introduction of the resonance method by K. Soundararajan in \cite{Sound08}, who showed that for the Riemann zeta function $\zeta(s)$ and sufficiently large $T$, 
\[
\max_{T\leq t\leq 2T}|\zeta( \frac 12+it)|\geq \exp\left((1+o(1))\sqrt{\frac{\log T}{\log_2 T}}\right),
\]
 where we denote $\log_j$ the $j$-fold iterated logarithm throughout the paper. The above result was improved by A. Bondarenko and K. Seip in \cite{BS18} using the long resonators introduced by C. Aistleitner in \cite{Aistleitner16} to be
\[
\max_{0\leq t\leq T}|\zeta(\frac 12+it)|\geq \exp\left((1+o(1))\sqrt{\frac{\log T \log_3 T}{\log_2 T}}\right). 
\] 
  The constant $1$ was further improved to $\sqrt{2}$ by R. de la Bret\`eche and G. Tenenbaum in \cite{BT19}.
  
 It is also shown in \cite[Theorem 2]{Sound08} that for sufficiently large $X$, 
\[
\max_{\substack{d \text{ fundamental discriminant }\\X<|d|\leq 2X}}|L(\frac 12, \chi_d)|\geq \exp\left(\left(\frac 1{\sqrt{5}}+o(1)\right)\sqrt{\frac{\log X }{\log_2 X}}\right),
\] 
 where $\chi_d=\leg {d}{\cdot}$ is the Kronecker symbol.

  In \cite{DM25}, P.  Darbar and G. Maiti employed the long resonator technique of Bondarenko--Seip to study the large value of $L(1/2, \chi_d)$ to show that  under the Generalized Riemann Hypothesis (GRH), for sufficiently large $X$, 
\begin{align}
\label{d}
	\max_{\substack{d \text{ fundamental discrimint }\\X<|d|\le 2X}}|L(\frac 12, \chi_{d})|\geq \exp\left(\left(\frac 1 2 +o(1) \right)\sqrt{\frac{\log X \log_3 X}{\log_2 X}}\right).
\end{align}

  Throughout the paper, we reserve the letters $p, q$ for primes.  In \cite{FHX26}, M. Fan, S. Hua and S. Xie  proved that for sufficiently large $X$, 
\begin{align}
\label{p}
\max_{\substack{X< p \leq 2X \\ p \equiv 1 \shortmod 8}}|L(\frac 12, \chi_p)|\geq \exp\left(\left(\sqrt{\frac 8{45}}+o(1)\right)\sqrt{\frac{\log X }{\log_2 X}}\right).
\end{align}

 It is the aim of this paper to carry the approach of Darbar and Maiti to study large values at the central point of the family of quadratic Dirichlet $L$-functions of prime moduli under GRH.  Due to some technical reasons, we consider the family of $L$-functions $\{ L(s, \chi_{8q}) \}$ with $q$ running over odd primes. Here we note as shown in \cite[Theorem 9.13]{MVa1}, the character $\chi_{8q}$ is primitive modulo $8q$. 
 
    Our main result investigates the maximum sizes of $ L(s, \chi_{8q})$ for $q$ running over primes.  
\begin{theorem}
\label{main theorem 1}
 With the notation as above and assuming the truth of GRH. We have for sufficiently large $X$, 
\begin{align}
\label{Llowerbound}
\begin{split}
  \max_{\substack{q \text{ \rm prime }\\X< q \le 2X}}\Big|L(\frac 12, \chi_{8q})\Big |\geq \exp\left(\left(\frac 1 2 +o(1) \right)\sqrt{\frac{\log X \log_3 X}{\log_2 X}}\right).
\end{split}
\end{align}	
\end{theorem}

   As mentioned above, our proof of Theorem \ref{main theorem 1} is based on the approaches introduced in \cite{DM25}. The result of Theorem \ref{main theorem 1} implies that given in \eqref{d}. Without much effort, one may show that Theorem \ref{main theorem 1} is also valid with $L(\frac 12, \chi_{8q})$ replaced by $L(\frac 12, \chi_{q})$ for $q$ being a fundamental discriminant so that our result also  improves that given in \eqref{p} under GRH.

\section{Preliminaries}
\label{sec 2}

\subsection{Approximate functional equation}
\label{sect: FE}
We have the following approximate functional equation concerning $L(1/2, \chi_{8q})$ given in \cite[Lemma 2.2]{sound1}.
\begin{lemma}
\label{lem:AFE}
  For any odd prime  $q$, we have
\begin{align*}
\begin{split}
 L\left( \frac{1}{2}, \chi_{8q} \right) = & 2\sum^{\infty}_{\substack{n=1}} \frac{\chi_{8q}(n)}{\sqrt{n}} V
\left(\frac{ n}{\sqrt{q}} \right),
\end{split}
\end{align*}
 where for any real number $x>0$,
\begin{align*}
 V(x) = \frac{1}{2 \pi i} \int\limits\limits_{(2)}  \left(\frac{8}{\pi}\right)^{s/2}
 \frac {\Gamma(s/2+1/4)}{\Gamma(1/4)}  x^{-s} \frac {\dif s}{s}.
\end{align*}
\end{lemma}

  As shown in \cite[Lemma 2.1]{sound1}, the function $V(x)$ is real-valued, smooth on $(0, +\infty)$, bounded as $x$ approaches $0$ and decays exponentially as $x\to +\infty$. More precisely,  we have for any $\varepsilon>0$, 
\begin{equation} 
\label{2.07}
      V\left (x \right) = 1+O(x^{1/2-\varepsilon}) \; \mbox{for} \; 0<x <1   \quad \mbox{and} \quad V^{(j)}\left (x \right) =O(e^{-x}) \; \mbox{for}
      \; x >0, \; j \geq 0.
\end{equation}

   We also note that
\begin{align}
\label{eq:Vder}
 V'(x) = -\frac{1}{2 \pi i} \int\limits\limits_{(2)}   \left(\frac{8}{\pi}\right)^{s/2}
 \frac {\Gamma(s/2+1/4)}{\Gamma(1/4)}  x^{-s-1} ds.
\end{align}   
   Recall that for $\Re(s)>0$,  
\begin{align*}
  \Gamma(s)=\int^{\infty}_0e^{-x}x^{s-1}dx.
\end{align*}      
  It follows from the inverse Mellin transformation that for $x>0$, 
\begin{align*}
  e^{-x}=\frac{1}{2 \pi i} \int\limits\limits_{(2)}\Gamma(s)x^{-s}ds. 
\end{align*}
  We deduce from the above that
\begin{align*}
  2x^{1/4}e^{-x}=\frac{1}{2 \pi i} \int\limits\limits_{(2)}\Gamma(\frac s2+\frac 14)x^{-s}ds. 
\end{align*} 
  It follows from  \eqref{eq:Vder} and the above that we have $V'(x)<0$ for $x>0$. Note moreover that by \eqref{2.07} we have $\lim_{x \rightarrow \infty}V(x)=0$. Thus we conclude that $V(x)$ is positive for $x \geq 0$. 

\subsection{Smoothed character sums}
\label{smoothsum}

   We define $\delta_{c=\square}=1$ if $c$ is a perfect square and $\delta_{c=\square}=0$ otherwise.  In the remainder of the paper, let $\Phi$ denote a smooth, non-negative function compactly supported on $[1,2]$ satisfying $\Phi(x) =1$ for $x\in [5/4,7/4]$. The Mellin transform of $\Phi(x)$ is wirtten as ${\widehat \Phi}(s)$ so that for any complex number $s$,
\begin{equation*}
{\widehat \Phi}(s) = \int\limits_{0}^{\infty} \Phi(x)x^{s}\frac {\dif x}{x}.
\end{equation*}
   We have the following result on the smoothed quadratic character sums.
\begin{lemma}
\label{lemma logd}
With the notation as above and assuming the truth of GRH. Let $c$ be a positive odd integer and $\Phi(X)$ be a smooth function fitting the above descriptions. Then for any $\varepsilon>0$,
\begin{equation*} 
 \sum_{(q,2)=1} (\log q) \chi_{8q}(c) \Phi \left( \frac {q}X \right) =  \delta_{c=\square}\widehat{\Phi}(1)X+O \left( X^{1/2+\varepsilon}\log
  (c+2) \right).
 \end{equation*}
\end{lemma}

\section{Proof of the Theorem \ref{main theorem 1}}
\subsection{Initial Treatments}

  Let $a\in (1, \infty), \delta\in (0,1)$ be fixed such that $1<a<\frac{1}{\delta}$, and let $N$ be a large number to be specified later. Let $[x]$ denote the largest integer not exceeding $x$ for any real $x$.  For $k=1,\ldots,[(\log_{2}N)^{\delta}]$, denote $\mathcal{P}_{k}$ the set of all primes $p$ such that $$e^{k}\log N \log_{2} N < p \le e^{k+1}\log N \log_{2} N.$$
  We further denote $ \mathcal{R}_{k}$ the set of positive integers  that have at least $\frac{a\log N}{k^{2}\log_{3}N}$ prime divisors in $\mathcal{P}_{k}$, and $ \mathcal{R}^{\prime}_{k}$ the set of positive integers   from $\mathcal{R}_{k}$ that have all their prime divisors in $\mathcal{P}_{k}$. 
  
  We let $\mathcal{P}$ be the union of $\mathcal{P}_{k}$ so that it is the set of all primes $p$ such that $$e\log N \log_{2} N< p \le e^{(\log_{2} N)^{\delta}} \log N \log_{2} N.$$ 
We define a multiplicative function $\psi$ supported on the set of square-free integers such that for any  $p \in\mathcal{P}$,
\begin{align}
\label{psidef}
\psi(p)=\sqrt{\frac{\log N \log_{2}N }{\log_{3}N}}p^{-1/2}\left(\log p-\log(\log N \log_{2}N)\right)^{-1}.
\end{align}
 We further define $\psi(p)=0$ for $p \notin \mathcal{P}$.
 
  We now set
\[
\mathcal{R}:=\rm{supp}(\psi)\setminus\displaystyle \bigcup_{k=1}^{[(\log_{2}N)^{\delta}]}\mathcal{R}_{k}.
\]  
  Then, we define the resonator for $L(\frac 1 2, \chi_{8q})$ for any prime $q$ to be the Dirichlet polynomial 
\[
R_q:= \sum_{m\in \mathcal{R}}\psi(m)\chi_{8q}(m).
\] 
 
  The following lemma, taken from \cite[Lemma 3]{DM25}, gives estimations on the sizes of $|\mathcal{R}|$, respectively.  
\begin{lemma}\label{Le1} 
	With the notation as above. Suppose that $1<a<\frac{1}{\delta}$. Then we have $|\mathcal{R}|\leq N$ for $N$ large enough.
\end{lemma}

  We now define 
\begin{align*}
\mathcal{A}_{N}:=\frac{1}{\displaystyle \sum_{m \in \mn } \psi(m)^{2} }\sum_{n\in \mn }\frac{\psi{(n)}}{\sqrt{n}}\sum_{\substack{l|n }}\psi(l)\sqrt{l}.
\end{align*}

   We end this section by including three results from \cite[ Lemma 1--3]{BS18} concerning various sums related to $\mathcal{A}_{N}$.
\begin{lemma}\label{Le2}
	With the notation as above.  We have as $N\to\infty$,
	\[\mathcal{A}_{N}\ge \exp\left((\delta + o(1))\sqrt{\frac{\log N\log_{3}N}{\log_{2}N}}\right).
	\]	
\end{lemma}

\begin{lemma}\label{Le3}
	With the notation as above.  We have as $N\to\infty$,
	\begin{align*}
	\frac{1}{\displaystyle\sum_{m\in \mn} \psi(m)^{2}}\sum_{\substack{n\in \mn \\  n\notin\mathcal{R}}}\frac{\psi{(n)}}{\sqrt{n}}\sum_{\substack{l|n }}\psi(l)\sqrt{l}=o(\mathcal{A}_{N}).
	\end{align*}	
	
\end{lemma}

\begin{lemma}\label{Le5}
	With the notation as above. We have as $N\to\infty$, for any $\varepsilon>0$, 
	\[ 
	\frac{1}{\displaystyle \sum_{m \in \mn} \psi(m)^{2}}\sum_{\substack{n\in \mathcal{R} }}\frac{\psi{(n)}}{\sqrt{n}}\sum_{\substack{l|n\\ l \le n/N^{\varepsilon}}}\psi(l)\sqrt{l}=o(\mathcal{A}_{N}),
	\]
	where the implicit constant only depends on $\varepsilon$.
\end{lemma}

\subsection{Completion of the proof}
	 Let $\Phi$ be the function described in Section \ref{smoothsum}.  We set
\begin{align*} 
 \mathcal{S}_{1}:=\sumstar_q (\log q)L(\frac{1}{2},\chi_{8q}) R_q^2 \Phi (\frac {q}{X}), \quad \mathcal{S}_{2}:=\sumstar_{q }(\log q)R_q^2\Phi (\frac {q}{X}).
\end{align*}
  As $R^2_d \geq 0$ and $L(\frac{1}{2},\chi_{8q}) \in \mr$, we then observe that
\begin{align}  
\label{maxlower}
 \max_{\substack{X< q \le 2X}}\Big|L(\frac 12, \chi_{8q})\Big |\ge \frac{\mathcal{S}_1}{\mathcal{S}_2}. 
\end{align}
It remains to establish a lower bound for $\mathcal{S}_{1}$ and an upper bound for $\mathcal{S}_{1}$. We first note that
\begin{align}
\label{S1eval}
\begin{split}
	\mathcal{S}_{1}=&\sumstar_q (\log q)L(\frac{1}{2},\chi_{8q}) R_q^2 \Phi (\frac {q}{X}) \\
	=&2\sum_{m, n\in \mathcal{R}} \psi(m) \psi(n)\sum_{l\geq 1 } \frac{1}{\sqrt{l}}\sumstar_{q}\chi_{8q}(lmn) V\left(\frac{l}{\sqrt{q}} \right).   
\end{split}  
\end{align}
 We now apply Lemma \ref{lemma logd} and partial summation to see that
\begin{align*}
\begin{split}
&	\sumstar_{q}\chi_{8q}(lmn) V\left(\frac{l}{\sqrt{q}} \right) \Phi (\frac {q}{X}) \\
=& \int^{2X}_X V\left(\frac{l}{\sqrt{t}} \right) d \left( t\delta_{lmn=\square}{\widehat \Phi}(1)+ O\left(t^{1/2+\varepsilon}\log
  (lmn+2) \right)  \right) \\
=& {\widehat \Phi}(1)X \delta_{lmn=\square}\int^{2}_1 V\left(\frac{l}{\sqrt{Xt}} \right) dt \\
&+V\left(\frac{l}{\sqrt{t}} \right)O\left(t^{1/2+\varepsilon}\log
  (lmn+2) \right)\Big |^{2X}_X -\int^{2X}_XO\left(t^{1/2+\varepsilon}\log
 (lmn+2) \right)V'\left(\frac{l}{\sqrt{t}}\right ) \frac {l}{2t^{3/2}}dt \\
=: & {\widehat \Phi}(1)X  \delta_{lmn=\square}\int^{2}_1 V\left(\frac{l}{\sqrt{Xt}} \right) dt+R. 
\end{split}
\end{align*}
 It follows from this and \eqref{S1eval} that 
\begin{align*}	
    \mathcal{S}_{1}={\widehat \Phi}(1)X
     \sum_{m, n \in \mathcal{R}}&\psi(m)\psi(n)\sum_{\substack{lmn =\square}} \frac{1}{\sqrt{l}}\int_{1}^{2}V\left(\frac{l}{\sqrt{Xt}} \right) dt+O\Bigg(\sum_{m, n\in \mathcal{R}} \psi(m) \psi(n)\sum_{l\geq 1 } \frac{1}{\sqrt{l}}R\Bigg).
\end{align*}
    In view of the rapid decay of $V$ and $V'$ given in \eqref{2.07}, we see that we may restrict the sum over $l$ to be $l \leq X^{1/2+\varepsilon}$ for any $\varepsilon>0$ in the error term above.  We further set $N= X^{\frac{1}{4}-5\varepsilon}$ for some $0<\varepsilon<1/20$ and observe that $|\mathcal{R}|\leq N$ from Lemma \ref{Le1}. It follows that
\begin{align*}	
    R \ll  X^{1/2+\varepsilon}.
\end{align*}    
   We deduce from the above that
\begin{align}	
\label{S1exp}
\begin{split}
    \mathcal{S}_{1}=& {\widehat \Phi}(1)X\sum_{m, n \in \mathcal{R}}\psi(m)\psi(n)\sum_{\substack{ lmn =\square}} \frac{1}{\sqrt{l}}\int_{1}^{2}V\left(\frac{N(l)}{\sqrt{Xt}} \right) dt+O\Bigg(X^{1/2 +\varepsilon}\sum_{\substack{l \leq X^{1/2+\varepsilon}}}\frac{1}{\sqrt{l}}  \left(\sum_{m\in \mathcal{R}}\psi(m)\right)^2\Bigg) \\
    =& {\widehat \Phi}(1)X\sum_{m, n \in \mathcal{R}}\psi(m)\psi(n)\sum_{\substack{ lmn =\square}} \frac{1}{\sqrt{l}}\int_{1}^{2}V\left(\frac{N(l)}{\sqrt{Xt}} \right) dt+O\Bigg(X^{3/4 +2\varepsilon}\left(\sum_{m\in \mathcal{R}}\psi(m)\right)^2\Bigg) \\
    =& {\widehat \Phi}(1)X\sum_{m, n \in \mathcal{R}}\psi(m)\psi(n)\sum_{\substack{ lmn =\square}} \frac{1}{\sqrt{l}}\int_{1}^{2}V\left(\frac{N(l)}{\sqrt{Xt}} \right) dt+O\Bigg(X^{3/4 +2\varepsilon}|\mathcal{R}|\sum_{m\in \mathcal{R}}\psi^2(m)\Bigg), 
\end{split}
\end{align}
   where the last estimation above follows from the Cauchy--Schwarz inequality. 
   
   As $\psi(n) \geq 0$ and  $V(x) \geq 0$, we may thus keep only the case $lm=n$ in the main term of \eqref{S1exp}. Together with the observation that $|\mathcal{R}|\leq N$ from Lemma \ref{Le1}, we see that
\begin{align}
\label{S1remainder} 
\mathcal{S}_{1} 
	 \ge & {\widehat \Phi}(1)X\sum_{n \in\mathcal{R}} \frac{\psi(n)}{\sqrt{n}}  \sum_{\substack{m \mid n }} \psi(m)\sqrt{m}\int_{1}^{2}V\left(\frac{n}{m\sqrt{Xt}} \right) dt+O\left(X^{3/4+2\varepsilon}N\left(\sum_{\substack{m\in \mathcal{R}}}\psi(m)^2\right)\right).   
\end{align}
  
 As $N= X^{\frac{1}{4}-5\varepsilon}$, it follows that when $m\ge n/N^{\varepsilon}$, we have
$0< \frac{n}{m\sqrt{Xt}}<1$ for $1 \leq t \leq 2$. We now apply the estimation $V(x)=1 + O\left(x^{\frac{1}{2}-\varepsilon}\right)$ given in \eqref{2.07} to deduce from the above that
\begin{align}
\label{S1lower}
\begin{split}
\mathcal{S}_{1} 
	 \ge & {\widehat \Phi}(1)X\sum_{n \in\mathcal{R}} \frac{\psi(n)}{\sqrt{n}} \sum_{\substack{m\mid n\\ m\ge n/N^{\varepsilon}}} \psi(m)\sqrt{m} \\ 
& +O\left(X^{3/4+2\varepsilon}\sum_{n \in\mathcal{R}} \psi(n)\sum_{\substack{m\mid n\\ m\ge n/N^{\varepsilon}}} \psi(m) \right)+O\left(X^{1-\varepsilon}\left(\sum_{\substack{m\in \mathcal{R}}}\psi(m)^{2}\right)\right).
\end{split}
\end{align}

   We observe from \eqref{psidef} that we have $0< \psi(p) \leq 1$ for $p \in  \mathcal{P}$. It follows from this that
\begin{align*}
\sum_{n \in\mathcal{R}} \psi(n)\sum_{\substack{m\mid n\\ m\ge n/N^{\varepsilon}}} \psi(m)  \ll \sum_{n \in\mathcal{R}} \psi(n)\sum_{\substack{ m\mid n}} \psi(m)=\sum_{m \in\mathcal{R}} \psi(m)\sum_{\substack{n \in\mathcal{R} \\ m\mid n}}\psi(n) \ll |\mathcal{R}|\sum_{m \in\mathcal{R}} \psi^2(m). 
\end{align*}  

  As  $|\mathcal{R}|\leq N$ from Lemma \ref{Le1} and $N= X^{\frac{1}{4}-5\varepsilon}$,  we deduce from the above and \eqref{S1lower} that
\begin{align*}
\begin{split}
\mathcal{S}_{1} 
	 \ge & {\widehat \Phi}(1)X\sum_{n \in\mathcal{R}} \frac{\psi(n)}{\sqrt{n}} \sum_{\substack{m\mid n\\ m\ge n/N^{\varepsilon}}} \psi(m)\sqrt{m}+O\left(X^{1-\varepsilon}\left(\sum_{\substack{m\in \mathcal{R}}}\psi(m)^{2}\right)\right).
\end{split}
\end{align*}

   We now apply above together with Lemma \ref{Le2}, Lemma \ref{Le3}, and Lemma \ref{Le5} to see that
\begin{align}
\label{S1lower2}
\begin{split}
	 \mathcal{S}_{1}\ge&
	\left(1+o(1)\right){\widehat \Phi}(1)X\exp\left((\delta+ o(1))\sqrt{\frac{\log N\log_{3}N}{\log_{2}N}}\right)\left(\sum_{\substack{m\in \mathcal{R}}}\psi(m)^{2}\right) +O\left(X^{1-\varepsilon}\left(\sum_{\substack{m\in \mathcal{R}}}\psi(m)^{2}\right)\right)\\
		 \ge  &\left(1+o(1)\right){\widehat \Phi}(1)X \exp\left(\left(\delta \sqrt{\frac 14-5\varepsilon}+o(1)\right)\sqrt{\frac{\log X\log_{3} X}{\log_{2} X}}\right) \left(\sum_{\substack{m\in \mathcal{R}}}\psi(m)^{2}\right).
\end{split}
\end{align}

Next, we proceed to establish an upper bound of $\mathcal{S}_2$. Again, applying Lemma \ref{lemma logd} and arguing similar to those that lead to \eqref{S1remainder} imply that
\begin{align*}
	\mathcal{S}_{2}=& \sumstar_{q }(\log q)R_q^2\Phi (\frac {q}{X})=\sum_{m, n \in \mathcal{R}}\psi(m) \psi(n) \sumstar_{q }(\log q)\chi_{8q}(m n)\Phi (\frac {q}{X})\\
	            =& {\widehat \Phi}(1)X \sum_{\substack{m, n \in \mathcal{R} \\ m n =\square}}\psi(m) \psi(n)
	             +O\Bigg(X^{\frac 1 2 + \varepsilon} \sum_{\substack{m,n \in \mathcal{R}}}\psi(m) \psi(n)\log
  (mn+2)\Bigg)\\
	            =& {\widehat \Phi}(1)X	\sum_{\substack{n \in \mathcal{R}}}\psi(n)^{2} + O\left(X^{1/2 +\varepsilon} N \sum_{\substack{m\in  \mathcal{R}}} \psi(m)^{2} \right),
\end{align*}
 where the last equality above follows by noting that $mn=\square$ implies that $m=n$ since $m$ and $n$ are square-free.	

Again using $|\mathcal{R}|\leq N= X^{\frac{1}{4}-5\varepsilon}$, we deduce from the above that 
\begin{align*}            
  \mathcal{S}_{2} \le \left(1+o(1)\right){\widehat \Phi}(1)X
	            \left(\sum_{\substack{m\in \mathcal{R}}}\psi(m)^{2}\right).
\end{align*}
Hence, for sufficiently large $X$ and arbitrary small $\varepsilon>0$, we deduce from \eqref{maxlower}, \eqref{S1lower2} and the above that 
 \begin{align*}
\max_{\substack{q \text{ \rm prime }\\X< q \le 2X}}\Big|L(\frac 12, \chi_{8q})\Big | &\geq \frac{\mathcal{S}_1}{\mathcal{S}_2}\geq \exp\left(\left(\delta \sqrt{\frac 14-5\varepsilon}+o(1)\right)\sqrt{\frac{\log X\log_{3} X}{\log_{2} X}}\right).
\end{align*}
 By letting $\delta \rightarrow 1^+$, we obtain the desired estimation given in  \eqref{Llowerbound}. This completes the proof of Theorem \ref{main theorem 1}.

\hspace{0.1in}

\noindent{\bf Acknowledgments.} P. G. is supported in part by NSFC grant 12471003.

\vspace*{.5cm}


\begin{thebibliography}{99}
\bib{Aistleitner16}{article}{
    AUTHOR = {Aistleitner, C.},
     TITLE = {Lower bounds for the maximum of the {R}iemann zeta function
              along vertical lines},
   JOURNAL = {Math. Ann.},
    VOLUME = {365},
      YEAR = {2016},
    NUMBER = {1-2},
     PAGES = {473--496},
}

\bib{BS18}{article}{
    AUTHOR = {Bondarenko, A.}, 
      AUTHOR = {Seip, K.}, 
     TITLE = {Extreme values of the {R}iemann zeta function and its
              argument},
   JOURNAL = {Math. Ann.},
    VOLUME = {372},
      YEAR = {2018},
    NUMBER = {3-4},
     PAGES = {999--1015},
}

\bib{DM25}{article}{
    AUTHOR = {Darbar, P.},
    AUTHOR = {Maiti, G.},
     TITLE = {Large values of quadratic {D}irichlet {$L$}-functions},
   JOURNAL = {Math. Ann.},
    VOLUME = {392},
      YEAR = {2025},
    NUMBER = {4},
     PAGES = {4573--4605},
}

\bib{BT19}{article}{
    AUTHOR = {de la Bret\`eche, R.},
    AUTHOR = {Tenenbaum, G.},
     TITLE = {Sommes de {G}\'al et applications},
   JOURNAL = {Proc. Lond. Math. Soc. (3)},
    VOLUME = {119},
      YEAR = {2019},
    NUMBER = {1},
     PAGES = {104--134},
}

\bib{FHX26}{article}{
     AUTHOR = {Fan, M.},
     AUTHOR = {Hua, S.},
     AUTHOR = {Xie, S.},
     TITLE = {Extreme central values of quadratic {D}irichlet
              {$L$}-functions with prime conductors},
   JOURNAL = {Q. J. Math.},
    VOLUME = {77},
      YEAR = {2026},
    NUMBER = {1},
     PAGES = {175--199},
}

\bib{MVa1}{book}{
  author = 	 {H. L. Montgomery},
   author = 	 {R. C. Vaughan},
  title = 	 {{M}ultiplicative number theory. {I}. {C}lassical theory},
  publisher = 	 {Cambridge University Press},
  year = 	 {2007},
  volume = 	 {97},
  series = 	 { Cambridge Studies in Advanced Mathematics},
  address = 	 {Cambridge}
}

\bib{sound1}{article}{
       author =   {K. Soundararajan},
  title =    {Nonvanishing of quadratic {D}irichlet {$L$}-functions at $s=\frac{1}{2}$},
  journal =      {Ann. of Math. (2)},
  year =     {2000},
  volume =   {152},
  number =   {2},
  pages =    {447-488}
}

\bib{Sound08}{article}{
      AUTHOR = {Soundararajan, K.},
     TITLE = {Extreme values of zeta and {$L$}-functions},
   JOURNAL = {Math. Ann.},
    VOLUME = {342},
      YEAR = {2008},
    NUMBER = {2},
     PAGES = {467--486},
}

\end{thebibliography}
\end{document}